\documentclass[11pt]{article}

\usepackage{amsmath}
\usepackage{latexsym}
\usepackage{amsfonts}
\usepackage{amssymb}
\usepackage{amscd}
\usepackage{oldgerm}		

\def\Cal{\mathcal}
\def\Bbb{\mathbb}

\define\M{\Cal M}
\define\Y{\Cal Y}

\define\V{\Cal V}

\define\E{\Cal E}

\define\W{\Cal W}
\define\F{\Cal F}
\define\D{\Cal D}
\define\calS{\Cal S}

\define\C{{\Bbb C}}
\define\R{{\Bbb R}}
\define\Q{{\Bbb Q}}
\define\Z{{\Bbb Z}}

\define\Gc{G_{\Bbb C}}
\define\gg{{\textswab{g}}}	
\define\half{\frac{1}{2}}
\define\hph{\hphantom}

\define\Im{\operatornamewithlimits{Im}}
\define\pd{\partial}
\define\rel{{}^r}

\define\a{\alpha}
\define\b{\beta}
\define\g{\gamma}
\define\l{\lambda}

\define\lam{\Lambda^{-1,-1}}

\newtheorem{remark}{Remark.\!\!} 
\newtheorem{nonumtheorem}{Theorem.\!\!} 
\newtheorem{nonumlemma}{Lemma.\!\!} 
\newtheorem{nonumcorollary}{Corollary.\!\!}

\newtheorem{theorem}{Theorem}[section]
\newtheorem{lemma}{Lemma}[section]
\newtheorem{definition}{Definition}[section]
\newtheorem{corollary}{Corollary}[section]
\newtheorem{example}{Example}[section]
\newtheorem{construction}{Construction}[section]

\newbox\qedbox
\setbox\qedbox=\hbox{$\Box$}

\newenvironment{proof}{\smallskip\noindent{\bf Proof.}\hskip \labelsep}
                        {\hfill\penalty10000\copy\qedbox\par\medskip}

\begin{document}

\begin{center}{\bf Degenerations of Mixed Hodge Structure}
\end{center}
\vskip .5cm

 \centerline{Gregory J Pearlstein}
\vskip 1cm

\section{Introduction}

\par In this paper we extend Schmid's Nilpotent Orbit Theorem to admissible
variations of graded-polarized mixed Hodge structure
$$
	\Cal V\to\Delta^*
$$
and derive analogs of the harmonic metric equations for variations of 
graded-polarized mixed Hodge structure.  The original motivation for the 
study of such variations rests upon the following observation \cite{Deligne3}:
Let $f:Z\to S$ be a surjective, quasi-projective morphism.  Then, the sheaf
$$
 	\Cal V = R^k_{f\ast}(\C)			
$$
restricts to a variation of graded-polarized mixed Hodge structure over some 
Zariski-dense open subset of $S$.  More recently, such variations have been
shown to arise in connection with the study of the monodromy representations 
of smooth projective varieties \cite{Hain} as well as certain aspects of 
mirror symmetry \cite{Deligne4}.
\vskip 10pt

\par The basic problem of identifying a good class of abstract variations of 
graded-polarized mixed Hodge structure for which one could expect to obtain 
analogs of Schmid's orbit theorems was posed by Deligne in \cite{Deligne3}.  
The accepted answer to this question was provided by 
Steenbrink and Zucker in \cite{SZ}, wherein they introduced the category of 
{\it admissible} variations of graded-polarized mixed Hodge 
structure, and proved that every geometric variation $\V$ defined over a
smooth, quasi-projective curve $X$ is admissible, and moreover the cohomology
$H^k(X,\V)$ of such a curve $X$ with coefficients in an admissible variation
$\V\to X$ carries a functorial mixed Hodge structure.  The original question 
of developing analogs of Schmid's orbit theorems for such variations has
however remained largely unresolved.
\vskip 10pt

\par The general outline of the paper is as follows: In \S 2, we review the 
basic properties of the classifying spaces of graded-polarized mixed Hodge 
structures
$$
	\M = \M(W,\calS,\{h^{p,q}\})
$$
constructed in \cite{higgs} and recall how the isomorphism class of
a variation of graded-polarized mixed Hodge structure $\V\to S$ may be 
recovered from the knowledge of its monodromy representation:
$$
	\rho:\pi_1(S,s_0)\to GL(\V_{s_0}),\qquad
	\text{Image}(\rho) = \Gamma
$$
and its period map
$$
	\varphi:S\to\M/\Gamma
$$
Following \cite{Kaplan}, we then construct a natural hermitian metric $h$ on 
$\M$ which is invariant under the action of the Lie group 
$$
	G_{\R} = \{\,g\in GL(V_{\R})^W \mid Gr(g)\in Aut(\calS,\R)\,\}
$$

\begin{remark} The Lie group $G_{\R}$ defined above only acts transitively 
upon the real points of $\M$ (i.e. the points $F\in\M$ for which the 
corresponding mixed Hodge structure $(F,W)$ is split over $\R$).    
\end{remark}

\par In \S 3, we derive analogs of the harmonic metric equations for filtered
vector bundles and determine necessary and sufficient conditions for such
a filtered harmonic metric to underlie a complex variation of graded-polarized
mixed Hodge structure.
\vskip 10pt

\par In \S 4, we recall the notion of an admissible variation of graded
polarized mixed Hodge structure and present a result of P. Deligne 
\cite{Deligne} which shows how to construct a distinguished 
$sl_2$-representation from the limiting data of an admissible 
variation $\V\to\Delta^*$.
\vskip 10pt

\par Making use of the material of \S 2 and \S 4, we then proceed in \S 5
to prove an analog of Schmid's Nilpotent Orbit Theorem for admissible 
variations of graded-polarized mixed Hodge structure $\V\to\Delta^*$ which 
gives distance estimates identical to those of the pure case.  Namely, once 
the period mapping of such a variation is lifted to a map of the upper 
half-plane into $\M$, the Hodge filtration of the given variation and that of 
the associated nilpotent orbit satisfy an estimate of the form 
$$
	d_{\M}(F(z),e^{zN}.F_{\infty})
	\leq K\Im(z)^{\b}\exp(-2\pi\Im(z))
$$

$$
	\text{\bf Acknowledgments}
$$

\par As much of this work stems from the author's thesis work, he would like
to thank both his advisor Aroldo Kaplan and his committee members David Cox
and Eduardo Cattani for their guidance.  The author would also like to thank 
Ivan Mirkovic and Richard Hain for enabling his stay at Duke University during
the 1998--1999 academic year, P. Deligne both for his contribution of the 
main lemma of \S 4 and his many helpful comments.

\setcounter{theorem}{0}
\setcounter{lemma}{0}
\setcounter{definition}{0}
\setcounter{corollary}{0}
\setcounter{equation}{0}
\setcounter{example}{0}

\section{Preliminary Remarks} In this section we review some background
material from \cite{SZ} and \cite{higgs}.

\setcounter{definition}{0}					
\begin{definition} Let $S$ be a complex manifold.  Then, a variation of 
graded-polarized mixed Hodge structure $\V\to S$ consists of a $\Q$-local 
system $\V_{\Q}$ over $S$ equipped with:
\begin{description}
\item{(1)} \hph{a}A rational, increasing weight filtration 
	   $0\subseteq\cdots\W_k\subseteq\W_{k+1}\subseteq\cdots
		 \subseteq\V_{\C}
	   $
	   of $\V_{\C} = \V_{\Q}\otimes\C$.
\item{(2)} \hph{a}A decreasing Hodge filtration 
	   $
		0\subseteq\cdots\F^p\subseteq\F^{p-1}\subseteq\cdots
		 \subseteq\V
	   $
           \linebreak of $\V = \V_{\C}\otimes{\Cal O}_S$ by holomorphic 
	   subbundles.
\item{(3)} \hph{a}A collection of rational, non-degenerate bilinear forms
	   $$
		\calS_k:Gr^{\W}_k(\V_{\Q})\otimes Gr^{\W}_k(\V_{\Q})\to\Q
	   $$
	   of alternating parity $(-1)^k$.
\end{description} 
satisfying the following mutual compatibility conditions:
\begin{description}
\item{(a)} \hph{a}Relative to the Gauss--Manin connection of $\V$:
	   $$
		\nabla\,\F^p\subseteq\Omega^1_S\otimes\F^{p-1}
	   $$
\item{(b)} \hph{a}For each index $k$, the triple 
	   $(Gr^{\W}_k(\V_{\Q}),\F\,Gr^{\W}_k(\V_{\Q}),\calS_k)$
	   defines a variation of pure, polarized Hodge structure
           of weight $k$.
\end{description}
\end{definition}

\par As discussed in \cite{higgs}, the data of such a variation $\V\to S$ may 
be effectively encoded into its monodromy representation
\setcounter{equation}{1}
\begin{equation}
        \rho:\pi_1(S,s_0)\to GL(\V_{s_0}),\qquad
	\text{Image}(\rho) = \Gamma			
\end{equation}
and its period map 
\begin{equation}
        \phi:S\to\M/\Gamma				
\end{equation}

\par To obtain such a reformulation, observe that it suffices to consider a
variation $\V$ defined over a simply connected base space $S$.  
Trivialization of $\V$ relative to a fixed reference fiber $V=\V_{s_0}$ via 
parallel translation will then determine the following data:
\begin{description}
\item{(1)} \hph{a}A rational structure $V_{\Q}$ on $V$.
\item{(2)} \hph{a}A rational weight filtration $W$ of $V$.
\item{(3)} \hph{a}A variable Hodge filtration $F(s)$ of $V$.
\item{(4)} \hph{a}A collection of rational, non-degenerate bilinear forms
            $$
		\Cal S_k:Gr^W_k\otimes Gr^W_k\to\C
	    $$ 
	    of alternating parity $(-1)^k$.
\end{description}
subject to the restrictions:  
\begin{description}
\item{(a)} \hph{a}The Hodge filtration $F(s)$ is holomorphic and horizontal, 
	   i.e.
           \begin{equation}
              \frac{\pd}{\pd\bar s_j}\,F^p(s)\subseteq F^p(s),\qquad
              \frac{\pd}{\pd s_j}\, F^p(s)\subseteq F^{p-1}(s) 
           \end{equation}
           relative to any choice of holomorphic coordinates $(s_1,\dots,s_n)$
           on $S$.
\item{(b)} \hph{a}Each pair $(F(s),W)$ is a mixed Hodge structure, 
	   graded-polarized by the bilinear forms $\{\Cal S_{k}\}$.
\end{description}

\par Conversely, the data listed in items $(1)$--$(4)$ together with the 
restrictions $(a)$ and $(b)$ suffice to determine a VGPMHS over a simply 
connected base.
\vskip 10pt

\par To extract from these properties an appropriate classifying space of 
graded-polarized mixed Hodge structures, observe that by virtue of conditions 
$(a)$ and $(b)$, the graded Hodge numbers $h^{p,q}$ of $\V$ are constant.  
Consequently, the filtration $F(s)$ must assume values in the set 
$$
	\M = \M(W,\Cal S,h^{p,q})
$$ 
consisting of all decreasing filtrations $F$ of $V$ such that
\begin{itemize}
\item $(F,W)$ is a mixed Hodge structure, graded-polarized by $\cal S$.
\item $\text{dim}_{\C}\,F^pGr^W_k=\sum_{r\geq p}\,h^{r,k-r}$.
\end{itemize}

\par To obtain a complex structure on $\M$, one simply exhibits $\M$ as an
open subset of an appropriate \lq\lq compact dual\rq\rq{} $\check\M$.  More
precisely, one starts with the flag variety $\check\F$ consisting of all 
decreasing filtrations $F$ of $V$ such that
$$
	\text{dim}\,F^p = f^p,\qquad f^p = \sum_{r\geq p,s}\,h^{r,s}
$$
To take account of the weight filtration $W$, one then defines $\check\F(W)$ 
to be the submanifold of $\check\F$ consisting of all filtrations 
$F\in\check\F$ which have the additional property that
$$
	\text{dim}\,F^p Gr^W_k = \sum_{r\geq p}\, h^{r,k-r}
$$
As in the pure case, the appropriate \lq\lq compact dual\rq\rq{} 
$\check\M\subseteq\check\F(W)$ is the submanifold of $\check\F(W)$ 
consisting of all filtrations $F\in\check\F(W)$ which satisfy Riemann's 
first bilinear relation with respect to the graded-polarizations $\calS$.
In particular, as shown in \cite{higgs}, $\check\M$ contains the 
classifying space $\M$ as a dense open subset.
\vskip 10pt

\par In order to state our next result, we recall that \cite{CKS} each 
choice of a mixed Hodge structure $(F,W)$ on a complex vector space 
$V=\V_{\R}\otimes\C$ determines a unique, functorial bigrading 
$$
	V = \bigoplus_{p,q}\, I^{p,q}
$$
with the following three properties:
\begin{description}
\item{$(1)$} For each index $p$, $F^p = \oplus_{a\geq p}\, I^{a,b}$.
\item{$(2)$} For each index $k$, $W_k = \oplus_{a+b\leq k}\, I^{a,b}$.
\item{$(3)$} For each bi-index $(p,q)$, 
	     $\bar I^{p,q} = I^{q,p} \mod \oplus_{r<q,s<p}\, I^{r,s}$.
\end{description}
\par In analogy with the pure case, I shall call this decomposition the 
Deligne--Hodge decomposition of $V$.  In particular, as discussed in 
\cite{higgs}, the pointwise application of this construction to a variation
of graded-polarized mixed Hodge structure $\V\to S$ determines a smooth
decomposition 
$$
	\V = \bigoplus_{p,q}\, {\Cal I}^{p,q}
$$
of $\V$ into a sum of $C^{\infty}$-subbundles.

\begin{nonumtheorem}[Kaplan, \cite{Kaplan}] Let $\M_{\R}$ denote the set of 
all filtrations $F\in\M$ for which the corresponding mixed Hodge structure 
$(F,W)$ is split over $\R$, (i.e. $\overline{I^{p,q}} = I^{q,p}$) and 
$G_{\R} = \{\, g\in GL(V_{\R})^W \mid Gr(g)\in Aut(\calS,\R) \,\}$.  Then, 
\begin{itemize}
\item The group $G_{\R}$ acts transitively on $\M_{\R}$.
\item The group $\Gc = \{\, g\in GL(V)^W \mid Gr(g)\in Aut(\calS,\C) \,\}$
      acts transitively on $\check\M$.
\item The intermediate group
      $G = \{\, g\in GL(V)^W \mid Gr(g)\in Aut(\calS,\R) \, \}$
      acts transitively on $\M$.
\end{itemize}
\end{nonumtheorem}

\begin{remark} By virtue of functoriality, each graded-polarized mixed Hodge
structure $(F,W)$ determines a natural mixed Hodge structure on the 
Lie group $\gg = Lie(\Gc)$ via the bigrading
$\gg^{r,s}_{(F,W)}
   = \{\, a\in\gg \mid \a:I^{p,q}_{(F,W)}\to I^{p+r,q+s}_{(F,W)} \,\}$.
\end{remark}

\par To measure distances in $\M$, we shall now describe the construction of
a natural hermitian metric $h$ carried by the classifying space $\M$.  For
the most part, our presentation follows the unpublished notes \cite{Kaplan}.
To motivate this construction, recall that in the pure case, there are 
essentially two equivalent ways of constructing a natural hermitian
metric on the classifying spaces $\D$:
\begin{description}
\item{(1)} \hph{a}Use the fact that $G_{\R}$ acts transitively on $\D$ with 
	    compact isotopy to define a $G_{\R}$-invariant metric on $\D$.
\item{(2)} \hph{a}Use the flag manifold structure of $\D$ and the Hodge metric 
	   $$
		h_F(u,v) = \calS(C_F u,\bar v)
	   $$
	   to induce a hermitian metric on $\D$.
\end{description}
The problem with using the first approach in the mixed case is that
\begin{itemize}
\item Although $G_{\R}^F$ is still compact, the group $G_{\R}$ does 
      in general act transitively upon $\M$.
\item Although $G$ does act transitively upon $\M$, this action 
      does not in general have compact isotopy.
\end{itemize}
Thus, in order to construct a natural hermitian metric on $\M$, we abandon 
the first approach and try instead to construct a natural generalization 
of the Hodge metric.  

\begin{nonumlemma}[Mixed Hodge Metric, \cite{Kaplan}] Let $(F,W,\Cal S)$ be 
a graded-polarized mixed Hodge structure with underlying vector space 
$V = V_{\R}\otimes\C$.  Then, there exists a unique, positive-definite 
hermitian inner product
$$
	h_F = h_{(F,W,\Cal S)}
$$
on $V$ with the following two properties:
\begin{description}
\item{$(a)$} The bigrading  $V=\bigoplus_{p,q}\,I^{p,q}_{(F,W)}$ is 
	     orthogonal with respect to $h_F$.
\item{$(b)$} If $u$ and $v$ are elements of $I^{p,q}$, then
	     $h_F(u,v) = i^{p-q}\Cal S_{p+q}([u],[\bar v])$.
\end{description}
\end{nonumlemma}

\setcounter{theorem}{4}					
\begin{theorem} The mixed Hodge metric $h$ defined above determines a natural 
hermitian metric on the classifying spaces of graded-polarized mixed Hodge 
structure 
$$
	\M = \M(W,\Cal S,h^{p,q})
$$
which is invariant under the action of $G_{\R}:\M\to\M$.
\end{theorem}
\begin{proof} Let $F$ be a point of the classifying space $\M$ and define
$$
	q_F = \bigoplus_{r<0,r+s\leq 0}\,\gg^{r,s}_{(F,W)} 
	\subseteq Lie(\Gc)
$$
Then, as discussed in \cite{higgs}, the subalgebra $q_F$ is a vector space
complement to $Lie(\Gc^F)$ in $Lie(\Gc)$, and hence the map
$\exp:u\in q_F \mapsto e^u.F\in\M$ restricts to a biholomorphism from 
some neighborhood of zero in $q_F$ onto some neighborhood of $F$ in $\M$.  
Consequently, we may introduce a hermitian metric on $T_F(\M)$ by first 
identifying $T_F(\M)$ with $q_F \cong T_0(q_F)$ via the differential
$\exp_*:T_0(q_F)\to T_F(\M)$ and then applying the rule:
\setcounter{equation}{5}
\begin{equation}
	h_F(\a,\b) = Tr(\a\b^*)					
\end{equation}
To see that $G_{\R}$ act by isometries with respect to $h$, one simply computes
using $(2.6)$.
\end{proof}

\setcounter{example}{6}						
\begin{example} Let $\M = \M(W,\Cal S,h^{p,q})$  be the classifying space of 
graded-polarized mixed Hodge structure defined by the following data
\begin{description}
\item{(1)} Rational structure: $V_{\Q} = \text{span}_{\Q}(e_0,e_2)$.
\item{(2)} Hodge numbers: $h^{1,1} = 1$, $h^{0,0} = 1$.
\item{(3)} Weight filtration: 
	   $$
		0 = W_{-1} \subset W_0 = \text{span}(e_0) 
		  = W_1	\subset W_2 = V
	   $$
\item{(4)} Graded Polarizations: $\Cal S_{2j}([e_j],[e_j]) = 1$.
\end{description}
Then, $\M$ is isomorphic to $\C$ via the map
$$
	\l\in\C \mapsto F(\l)\in\M,\qquad F^1(\l) = \text{span}(e_2 + \l e_0)
$$
Moreover, relative to the mixed Hodge metric, the following frame is both
holomorphic and unitary:
$$
	v_1(\l) = e_0,\qquad v_2(\l) = e_2 + \l e_0
$$
Consequently, the classifying space $\M$ is necessarily flat (relative to
the mixed Hodge metric).
\end{example}

\begin{remark} This example also shows that (in contrast with the pure case)
the period map $\phi:\Delta^*\to\M/\Gamma$ of an abstract variation of 
graded-polarized mixed Hodge structure $\V\to\Delta^*$ may have an irregular 
singularity at $s=0$.  To wit: Set 
	   $\Gamma = \{1\}$ and define
	   \setcounter{equation}{7} 
	   \begin{equation}
		\phi(s) = \exp(\frac{1}{s}):\Delta^*\to\M	
	   \end{equation}
	   where $\M\cong\C$ is the classifying space of Example $(2.7)$.
\end{remark}


\par As a prelude to \S 4, we now recall some background material from
\cite{SZ} and record some auxiliary results by the relative weight
filtration $\rel W$ and its relationship with the finite dimensional 
representations of $sl_2(\C)$.  To this end, we shall let $V$ denote 
a finite dimensional $\C$-vector space and define $E_{\a}(Y)$ to be 
the $\a$-eigenspace of a semisimple endomorphism $Y:V\to V$.

\setcounter{theorem}{8}					    
\begin{theorem} Given a nilpotent endomorphism $N:V\to V$, there exists a 
unique {\rm monodromy weight filtration} 
$$
	0\subset   W(N)_{-k} \subseteq W(N)_{1-k} \subseteq \cdots
	 \subseteq W(N)_{k-1} \subseteq W(N)_k = V
$$
of $V$ such that:
\begin{itemize}
\item $N:W(N)_j\to W(N)_{j-2}$ for each index $j$.
\item The induced maps 
      $N^j:Gr^{W(N)}_{j}\to Gr^{W(N)}_{-j}$ are isomorphisms.
\end{itemize}
\end{theorem}

\setcounter{example}{9}					   
\begin{example} Let $\rho$ be a finite dimensional representation of 
$sl_2(\C)$ and
$$
	N_{\pm} = \rho(n_{\pm}),\qquad Y = \rho(y)
$$
denote the images of the standard generators $(n_-,y,n_+)$ of $sl_2(\C)$.
Then, by virtue of the semisimplicity of $sl_2(\C)$ and the commutator
relations
$$
	[Y,N_{\pm}] = \pm 2N_{\pm}, \qquad [N_+,N_-] = Y
$$
it follows that: 
\setcounter{equation}{10}
\begin{equation}
	W(N_-)_k = \bigoplus_{j\leq k}\, E_j(Y)			
\end{equation}
\end{example}

\par To obtain a converse of the construction given in the preceding example,
recall that a grading $Y$ of an increasing filtration $W$ of a finite
dimensional vector space $V$ may be viewed as a semisimple element of $End(V)$ 
such that 
$$
	W_k = E_k(Y)\oplus W_{k-1}
$$
for each index $k\in\Z$.  In particular, each mixed Hodge structure $(F,W)$
defines a functorial grading $Y_{(F,W)}$ of the underlying weight filtration
$W$ via the rule
$$
	v\in I^{p,q}_{(F,W)} \implies Y_{(F,W)}(v) = (p+q)v
$$
Moreover, given any increasing filtration $W$ of a complex vector space $V$,
the set $\Y(W)\subset End(V)$ consisting of all gradings $Y$ of $W$ is an 
affine space upon which the nilpotent Lie algebra
$$
	Lie_{-1} = \{\, \a\in End(V) 
			\mid \a:W_k\to W_{k-1}\hph{a} \forall\,k\,\}
$$
acts transitively \cite{CKS}.

\setcounter{theorem}{11}				    
\begin{theorem} Let $N$ be a fixed, non-trivial, nilpotent endomorphism of
$V$, and 
$$
	(n_-,y,n_+)
$$ 
denote the standard generators of $sl_2(\C)$.  Then, there exists a bijective
correspondence between:
\begin{description}
\item{(a)} \hph{a}The set $S$ of all gradings $H$ of the monodromy 
	   weight filtration $W(N)$ for which 
	   $$
		[H,N] = -2N
	   $$
\item{(b)} \hph{a}The set $S^{\prime}$ of all representations
	   $$
		\rho:sl_2(\C) \to End(V)
	   $$
	   such that $\rho(n_-) = N$.
\end{description}
\end{theorem}
\begin{proof} The construction of Example $(2.10)$ defines a map
$f:S^{\prime}\to S$ by virtue of $(2.11)$ and the standard commutator relations
of $sl_2(\C)$.  To obtain an map $g:S\to S^{\prime}$ such that 
$$
	f\circ g = Id,\qquad g\circ f = Id
$$
one simply considers how the Jordan blocks of $N$ interact with the grading 
$H$ of $W(N)$.
\end{proof}

\setcounter{definition}{12}					
\begin{definition} Given an increasing filtration $W$ of $V$ and an integer
$\ell\in\Z$ the corresponding shifted object $W[\ell]$ is the increasing 
filtration of $V$ defined by the rule:
$$
		W[\ell]_j = W_{j+\ell}
$$
\end{definition}

\setcounter{theorem}{13}				     
\begin{theorem} Let $W$ be an increasing filtration of $V$.  Then, given a 
nilpotent endomorphism $N:V\to V$ which preserves $W$, there exists at most 
one increasing filtration of $V$ 
$$
	\rel W = \rel W(N,W)
$$
with the following two properties:
\begin{itemize}
\item For each index $j$, $N:\rel W_j\to\rel W_{j-2}$
\item For each index $k$, $\rel W$ induces on $Gr^W_k$ the corresponding 
shifted monodromy weight filtration
$$                  
			W(N:Gr^W_k\to Gr^W_k)[-k]
$$
\end{itemize}
\end{theorem}

\par Following \cite{SZ}, we shall call $\rel W(N,W)$ the relative
weight filtration of $W$ and $N$.  To relate this filtration with the
finite dimensional representations of $sl_2(\C)$, we may proceed as 
follows:

\setcounter{theorem}{14}				     
\begin{theorem} Let $\rel Y$ be a grading of $\rel W = \rel W(N,W)$ such that
\begin{itemize}
\item $\rel Y$ preserves $W$.
\item $[\rel Y,N] = -2N$.
\end{itemize}
Then, each choice of grading $Y$ of $W$ which preserves $\rel W$ determines
a corresponding representation $\rho:sl_2(\C)\to End(V)$ such that
$$
	\rho(n_-) = N_0,\qquad \rho(y) = \rel Y-Y
$$
where 
\setcounter{equation}{15}
\begin{equation}
	N = N_0 + N_{-1} + N_{-2} + \cdots		       
\end{equation}
denotes the decomposition of $N$ relative to the eigenvalues of $ad\,Y$. 
\end{theorem}
\begin{proof} By virtue of the definition of the relative weight filtration 
and the mutual compatibility of $\rel Y$ with $Y$, it follows that the induced
map
$$
	H = Gr(\rel Y - Y) : Gr^W \to Gr^W
$$
grades the monodromy weight filtration $W(N:Gr^W\to Gr^W)$.  Moreover, by 
hypothesis
$$
	[H,Gr(N)] = -2 Gr(N)
$$
Thus, application of Theorem $(2.12)$ defines a collection of representations
$$
	\rho_k:sl_2(\C)\to Gr^W_k
$$
which we may then lift to the desired representation $\rho:sl_2(\C)\to End(V)$ 
via the grading $Y$.

\par To see that $\rho(n_-) = N_0$ and $\rho(y) = \rel Y - Y$, observe that
if
$$
	\a = \a_0 + \a_{-1} + \a_{-2} + \cdots
$$
is the decomposition of an endomorphism $\a\in End(V)^W$  according to the 
eigenvalues of $ad\,Y$ then the lift of 
$$
	Gr(\a):Gr^W\to Gr^W
$$
with respect to the induced isomorphism $Y:Gr^W \to V$ is exactly $\a_0$.
Thus $\rho(n_-) = N_0$.  Likewise, $\rho(y) = \rel Y - Y$ since the mutual
compatibility of the gradings $\rel Y$ and $Y$ implies that $\rel Y$ and
$Y$ may be simultaneously diagonalized.
\end{proof}

\par To close this section, we now recall the following definition
from \cite{SZ}:

\setcounter{definition}{16}					
\begin{definition} A variation $\V\to\Delta^*$ of graded-polarized mixed 
Hodge structure with unipotent monodromy is said to be admissible provided:
\begin{description}
\item{$(i)$}  \hph{a}The limiting Hodge filtration $F_{\infty}$ of $\V$ exists.
\item{$(ii)$} \hph{a}The relative weight filtration $\rel W = \rel W(N,W)$ of 
the monodromy logarithm $N$ and the weight filtration $W$ of $\V$ exists.
\end{description}
\end{definition}

\par The precise meaning of condition $(i)$ is a follows:  Let
$$
	\varphi:\Delta^*\to\M/\Gamma
$$ 
be the period map of a variation of graded-polarized mixed Hodge structure
with unipotent monodromy.  Then, as in the pure case, $\varphi(s)$ may be 
lifted to a holomorphic, horizontal map
$$
	F(z):U\to\M
$$
from the upper half-plane $U$ into $\M$ which makes the following diagram
commute:
$$
\begin{CD}
                  U        @> F >>    \M        \\
	@Vs=\exp(2\pi i z)VV                    @VVV     \\
                   \Delta^* @> \phi >> \M/\Gamma
\end{CD}
$$
In particular, on account of this diagram, 
$$
	F(z+1) = e^N.F(z)
$$
and hence $F(z)$ descends to a \lq\lq untwisted period map\rq\rq{}
$$
	\psi(s):\Delta^*\to\check\M
$$
In accord with \cite{Schmid}, the limiting value of $\psi(s)$ at zero
[when it exists] is then called the limiting Hodge filtration $F_{\infty}$ 
of $\V$.

\begin{remark} To be coordinate free, the limiting Hodge filtration 
$F_{\infty}$ constructed above should be viewed as an object attached
to $T_0(\Delta)$.
\end{remark}

\setcounter{theorem}{0}
\setcounter{lemma}{0}
\setcounter{definition}{0}
\setcounter{corollary}{0}
\setcounter{equation}{0}
\setcounter{example}{0}

\section{Higgs Fields and Harmonic Metrics}

\par In this section we derive analogs of the harmonic metric equations for 
filtered vector bundles and discuss their relationship with complex 
variations of graded-polarized mixed Hodge structure.  

\begin{definition} A Higgs bundle consists of a $C^{\infty}$ complex vector
bundle endowed with a smooth linear map $\theta:E\to E\otimes\E^1$ of type
$(1,0)$ and a smooth differential operator $\bar\pd$ of type $(0,1)$ which 
satisfy the following mutual compatibility condition:
$$
	(\bar\pd + \theta)^2 = 0
$$
\end{definition}

\par Equivalently, by virtue of the Newlander--Nirenberg theorem, a Higgs 
bundle may be viewed as holomorphic vector bundle $(E,\bar\pd)$ endowed 
with the choice of a holomorphic map $\theta:E\to E\otimes\Omega^1$ such that 
$\theta\wedge\theta = 0$.

\setcounter{construction}{1}				
\begin{construction} Let $(E,\nabla)$ be flat vector bundle.  Then, each
choice of a hermitian metric $h$ determines a unique choice of differential
operators $\delta^{\prime}$ and $\delta^{\prime\prime}$ of respective types 
$(1,0)$ and $(0,1)$ such that the corresponding connections
$\delta^{\prime} + \nabla^{\prime\prime}$ and 
$\nabla^{\prime} + \delta^{\prime\prime}$ preserve the metric $h$.  In 
particular, each choice of a hermitian metric $h$ on a flat vector bundle 
$E$ determines a canonical decomposition
$$
	\nabla = \bar\theta + \bar\pd + \pd + \theta
$$
of the underlying flat connection $\nabla$ via the rule:
$$
   \bar\theta = \half(\nabla^{\prime\prime} - \delta^{\prime\prime}),\quad
   \bar\pd = \half(\nabla^{\prime\prime} + \delta^{\prime\prime}),\quad
   \pd = \half(\nabla^{\prime} + \delta^{\prime}),\quad
   \theta = \half(\nabla^{\prime}-\delta^{\prime})
$$
\end{construction}

\par Motivated by \cite{Simpson}, we define a hermitian metric $h$ on a flat 
vector bundle $E$ to be {\it harmonic} provided the resulting operator 
$\bar\pd + \theta$ produced by Construction $(3.2)$ is of Higgs--type
[i.e. $(\bar\pd+\theta)^2 = 0$ ].  

\setcounter{example}{2}						
\begin{example} The Hodge metric $h$ carried by a variation of pure, polarized
Hodge structure $\V\to S$ is harmonic.
\end{example}

\begin{nonumlemma} Let $X$ be a compact K\"ahler manifold.  Then, a flat
vector bundle $E\to X$ admits a harmonic metric $h$ if and only if it is
semisimple.  
\end{nonumlemma}

\par In particular, since monodromy representation of a variation of
graded-polarized mixed Hodge structure need not be semisimple, 
the underlying flat vector bundle of a VGPMHS is need not admit a 
harmonic metric.  Accordingly, we devote the remainder of this section 
to deriving  analogs of the harmonic metric equations for flat, filtered 
vector bundles.

\setcounter{definition}{3}					
\begin{definition} A filtered bundle $(E,\W)$ consists of a $C^{\infty}$ 
vector bundle $E$ endowed with an increasing filtration
$$
   0 \subseteq\cdots\subseteq\W_{k-1}\subseteq\W_k\subseteq\cdots\subseteq E
$$
of $E$ by $C^{\infty}$ subbundles.  Likewise, a flat, filtered bundle 
$(E,\W,\nabla)$ consists of a $C^{\infty}$ vector bundle $E$ endowed with 
an increasing filtration by flat subbundles $\W_k\subseteq E$.
\end{definition}

\setcounter{construction}{4}			
\begin{construction} Each choice of a hermitian metric $h$ on a filtered bundle
$(E,\W)$ defines a corresponding grading $\Y_h$ of the underlying weight
filtration $\W$ by simply declaring $E_k(\Y_h)$ to be the orthogonal 
complement of $\W_{k-1}$ in $\W_k$ with respect to $h$.
\end{construction}

\par In particular, each choice of a hermitian metric $h$ on a flat, filtered
bundle $(E,\W,\nabla)$ determines a unique decomposition of the underlying
flat connection $\nabla$ into a sum of three terms 
\setcounter{equation}{5}
\begin{equation}
	\nabla = d_h + \tau_- + \theta_-			
\end{equation}
such that 
\begin{itemize}
\item $d_h$ is a connection which preserves the grading $\Y_h$.
\item $\tau_-$ is a nilpotent tensor field of type $(0,1)$ which maps
      $\W_k$ to $\W_{k-1}$ for each index $k$.
\item $\theta_-$ is a nilpotent tensor field of type $(1,0)$ which maps
      $\W_k$ to $\W_{k-1}$ for each index $k$.
\end{itemize}
Moreover, the resulting connection $d_h$ is automatically flat, as may be 
seen by expanding out the flatness condition $\nabla^2 = 0$ and taking note 
of the fact that $d_h:E_k(\Y_h)\to E_k(\Y_h)\otimes\E^1$ whereas
$\tau_-$, $\theta_-:\E_k(\Y_h)\to \W_{k-1}\otimes\E^1$.
\vskip 10pt

\par To continue, we note that
\begin{itemize}
\item To each filtered bundle $(E,\W)$ one may associate a corresponding graded
vector bundle $Gr^{\W}$ via the rule:
$$
	Gr^{\W}_k = \frac{\W_k}{\W_{k-1}} 
$$
\item By virtue of the induced isomorphism $\Y_h:E\cong Gr^{\W}$, the choice 
of a hermitian metric $h$ on a filtered bundle $(E,\W)$ determines an 
corresponding hermitian metric $h_{Gr^{\W}}$ on $Gr^{\W}$.
\end{itemize}

\setcounter{definition}{6}				
\begin{definition} A hermitian metric $h$ on a flat, filtered bundle 
$(E,\W,\nabla)$ is said to be graded-harmonic provided it is harmonic
with respect to the flat connection $d_h$ (i.e. the induced metric 
$h_{Gr^{\W}}$ on $Gr^{\W}$ is harmonic).
\end{definition}

\par In particular, the choice of a graded-harmonic metric $h$ on a flat 
filtered bundle $(E,\W,\nabla)$ determines an associated decomposition 
\setcounter{equation}{7}
\begin{equation}
	d_h = \bar\theta_0 + \bar\pd_0 + \pd_0 + \theta_0	
\end{equation}
by application of Construction $(3.2)$ to the pair $(h,d_h)$.
 
\setcounter{definition}{8}					
\begin{definition} Let $(E,\W,\nabla)$ be a flat, filtered bundle.  Then, a
hermitian metric $h$ on $E$ is said to be filtered harmonic provided:
\begin{description}
\item{$(a)$} $h$ is graded-harmonic.
\item{$(b)$} The differential operator $\bar\pd + \theta$ on $E$ obtained
by setting
$$
	\bar\pd = \bar\pd_0 + \tau_-,\qquad \theta = \theta_0 + \theta_-
$$
defined by equations $(3.6)$ and $(3.8)$ is of Higgs--type.
\end{description}
\end{definition}

\par To relate the preceding constructions to variations of graded-polarized
mixed Hodge structure, we recall from \S 2 that such a variation $\V$ comes
equipped with a canonical $C^{\infty}$ bigrading
\setcounter{equation}{9}
\begin{equation}
	\V = \bigoplus_{p,q}\, {\Cal I}^{p,q}			
\end{equation}
such that:
\begin{itemize}
\item For each index $p$, $\F^p = \oplus_{a\geq p}\, {\Cal I}^{a,b}$.
\item For each index $k$, $\W_k = \oplus_{a+b\leq k}\, {\Cal I}^{a,b}$.
\item For each bi-index $(p,q)$, 
$\overline{{\Cal I}^{p,q}} 
  = {\Cal I}^{q,p} \mod \bigoplus_{r<q,s<p}\, {\Cal I}^{r,s}$.
\end{itemize}

\par In particular, on account of this bigrading, a variation of 
graded-polarized mixed Hodge structure $\V$ comes equipped with 
the following additional structures:
\begin{itemize}
\item A canonical grading $\Y$ of $\W$ which acts as multiplication by $p+q$
on ${\Cal I}^{p,q}$.
\item A canonical mixed Hodge metric $h$ defined by pulling back the Hodge
metric of $Gr^{\W}$ via the induced isomorphism $\Y:\V\cong Gr^{\W}$.
\item A canonical Higgs bundle structure $\bar\pd + \theta$ (cf. \cite{higgs}).
\end{itemize}

\par To explain the construction of the Higgs bundle structure 
$\bar\pd + \theta$ on the underlying $C^{\infty}$ vector bundle of $\V$, 
we recall that a complex variation of Hodge structure ($\C$VHS) consists of 
a flat vector bundle $(E,\nabla)$ equipped with a $C^{\infty}$-decomposition
\begin{equation}
	E = \bigoplus_p\, E^p					
\end{equation}
which satisfies the horizontality condition
\begin{equation}
  \nabla:\E^0(E^p) \to 
  \E^{0,1}(E^{p+1})\oplus\E^{0,1}(E^p)\oplus\E^{1,0}(E^p)
	\oplus\E^{1,0}(E^{p-1})
\end{equation}
and note that a parallel hermitian bilinear form $S$ is then said to polarize
such a variation $E$ provided:
\begin{itemize}
\item The direct sum decomposition $(3.11)$ is orthogonal with respect to $S$.
\item The associated \lq\lq Hodge metric\rq\rq{}   
      $$
		h(u,v) = (-1)^p S(u,v),\qquad u,v \in E^p
      $$
      is positive definite.
\end{itemize}
Moreover, as with variations of pure, polarized Hodge structure, the Hodge 
metric of a polarized $\C$VHS is harmonic. 

\setcounter{lemma}{12}						
\begin{lemma} Let $E$ is a complex variation of Hodge structure and
$\bar\pd+\theta$ be the differential operator on $E$ obtained by 
decomposing the underlying flat connection
$$
	\nabla = \bar\theta + \bar\pd + \pd + \theta
$$
of $E$ in accord with equation $(3.12)$.  Then, $(\bar\pd + \theta)^2 = 0$ and
hence $\bar\pd + \theta$ is an operator of Higgs--type.
\end{lemma}
\begin{proof} One simply expands out the flatness condition $\nabla^2=0$ and
taking note of the additional requirements imposed by equation $(3.11)$.
\end{proof}

\setcounter{theorem}{13}				
\begin{theorem} Let $\V$ be a variation of graded-polarized mixed Hodge
structure.  Then, the $C^{\infty}$ decomposition
$$
	\V = \bigoplus_p\, E^p,\qquad E^p = \bigoplus_q\, {\Cal I}^{p,q}
$$
defines a complex variation of (unpolarized) Hodge structure on the underlying
flat bundle of $\V$.  
\end{theorem}

\begin{nonumcorollary} A variation of graded-polarized mixed Hodge structure
carries a canonical Higgs bundle structure $\bar\pd+\theta$.
\end{nonumcorollary}

\par In order to obtain a partial converse of the Lemma $(3.13)$, recall that 
a decomposition
$$
	E = \bigoplus_p\, E^p
$$
of a Higgs bundle $(E,\bar\pd + \theta)$ into a sum of holomorphic subbundles
is said to be a system of Hodge bundles if and only if 
$$
	\theta:E^p\to E^{p-1}\otimes\Omega^1
$$ 

\setcounter{lemma}{14}					
\begin{lemma} A Higgs bundle $(E,\bar\pd + \theta)$ defined over a compact
complex manifold $X$ admits a decomposition into a system of Hodge bundles 
if and only if 
$$
	(E,\bar\pd+\theta) \cong (E,\bar\pd + \l\theta)
$$
for each element $\l\in\C^*$.
\end{lemma}
\begin{proof} If $(E,\bar\pd+\theta)$ admits a decomposition into a sum of 
Hodge bundles $E = \bigoplus_p\, E^p$ then the desired isomorphism 
$f:(E,\bar\pd + \theta)\cong(E,\bar\pd + \l\theta)$ may be obtained by
setting $f=\l^p$ on $E^p$.

\par Conversely, given a Higgs bundle $(E,\bar\pd+\theta)$ which is a fixed
point of the $\C^*$ action $(E,\bar\pd + \theta) \to (E,\bar\pd + \l\theta)$
one may obtained the desired decomposition of $E$ into a system of Hodge
bundles as follows:  Let $f$ be an isomorphism from $(E,\bar\pd+\theta)$
to $(E,\bar\pd+\l\theta)$ for some element $\l\in\C^*$ which is not a
root of unity.  Then, because $f$ is holomorphic and $X$ is compact, 
the characteristic polynomial of $f$ is constant.  Upon decomposing 
$E$ into a sum of generalized eigenspaces, one then obtains the desired
system of Hodge bundles.
\end{proof}

\setcounter{corollary}{15}			
\begin{corollary}[Lemma 4.2, \cite{Simpson}] Let $X$ be a compact K\"ahler 
manifold.  Then, a semisimple flat bundle $E\to X$ underlies a polarized 
$\C$VHS if and only if the corresponding Higgs bundle $(E,\bar\pd + \theta)$ 
is a fixed point of the $\C^*$-action 
$$
	(E,\bar\pd + \theta) \mapsto (E,\bar\pd + \l\theta)
$$
Moreover, there is a 1-1 correspondence between the set of possible 
complex variations of polarized Hodge structure on such a semisimple
flat bundle $E$ and the set of possible decompositions of 
$(E,\bar\pd+\theta)$ into a system of Hodge bundles.
\end{corollary}

\par Motivated by the computations of \cite{higgs}, we now make the following
definition:

\setcounter{definition}{16}				
\begin{definition} A complex variation of mixed Hodge structure consists of
a flat vector bundle $(E,\nabla)$ endowed with a smooth decomposition of 
the underlying $C^{\infty}$-vector bundle into a sum of subbundles
\setcounter{equation}{17}
$$
	E = \bigoplus_{p,q}\, E^{p,q}			
$$
such that
\begin{description}
\item{$(a)$} For each index $k$, the sum 
$\W_k = \bigoplus_{p+q\leq k}\, E^{p,q}$ is a flat subbundle of $E$.
\item{$(b)$} If $E^p = \bigoplus_q\, E^{p,q}$ then for each bi-index $(p,q)$,
$$
	\nabla:\E^0(E^{p,k-p})\to\E^{0,1}(E^{p+1,k-p+1})\oplus
	\E^{0,1}(E^p)\oplus\E^{1,0}(E^{p,k-p})\oplus\E^{1,0}(E^{p-1})
$$
\end{description}
\end{definition}

\begin{remark} It is perhaps more natural to define a $\C$VMHS as consisting
of a $\C$-local system $\V$ endowed with a triple of suitably opposed 
filtrations $(\F,\W,\bar\F)$ which satisfy the requisite horizontality
conditions.  One then obtains a functor from this category to the category
of $\C$VMHS defined above via the map 
$$
	(\F,\W,\bar\F) \mapsto \bigoplus_{p,q}\, {\Cal I}^{p.q}
$$
However, in the absence of additional data (such as a real structure for $\V$),
one can only recover $\F$ and $\W$ from the ${\Cal I}^{p,q}$'s, and not the
full triple $(\F,\W,\bar\F)$.
\end{remark}

\par To define the notion of a complex variation of graded-polarized mixed
Hodge structure, note that by virtue of conditions $(a)$ and $(b)$, a complex 
variation of mixed Hodge structure induces complex variations of pure Hodge 
structure on $Gr^{\W}$ via the rule
$$
	E^p Gr^{\W}_k = \frac{E^p\cap\W_k + \W_{k-1}}{\W_{k-1}}
$$

\setcounter{definition}{17}					
\begin{definition} A complex variation of graded-polarized mixed Hodge
structure consists of a complex variation of mixed Hodge structure together
with a collection of parallel hermitian bilinear forms on $Gr^{\W}$ which
polarize the induced variations.
\end{definition}

\par In particular, every complex variation of graded-polarized mixed Hodge
structure carries a canonical mixed Hodge metric obtained by simply pulling 
back the Hodge metrics of $Gr^{\W}$ via the grading 
\setcounter{equation}{18}
\begin{equation}
	E_k(\Y) = \bigoplus_{p+q=k}\, E^{p,q}			
\end{equation}

\setcounter{theorem}{19}					
\begin{theorem} The Deligne--Hodge decomposition 
$$
	\V = \bigoplus_{p,q}\, {\Cal I}^{p,q}
$$
of a variation of graded-polarized mixed Hodge structure $\V\to S$ defines a 
complex variation of graded-polarized mixed Hodge structure with respect
to the underlying flat structure of $\V$.
\end{theorem}
\begin{proof} Since the weight filtration $\W$ of $\V$ is by definition flat, 
it will suffice to verify condition $(b)$ of Definition $(3.17)$ by simply 
computing the action of the Gauss--Manin connection $\nabla$ of $\V$ upon a 
smooth section $\sigma$ of ${\Cal I}^{p,q}$ at an arbitrary point 
$s_0\in S$. 

\par Accordingly, we recall from \S 2 that $\V$ may be represented near
$s_0$ by a holomorphic, horizontal map
$$
	F:\Delta^n\to\M
$$
from the polydisk $\Delta^n$ into a suitable classifying space $\M$ upon
selection a system of holomorphic local coordinates $(s_1,\dots,s_n)$ on
$S$ which vanish at $s_0$.  In particular, if $F=F(0)$ and 
$$
	q_F = \bigoplus_{r<0,r+s\leq 0}\, \gg^{r,s}
$$ 
denotes the vector space complement to $Lie(\Gc^F)$ in $Lie(\Gc)$ constructed
in \S 2 [cf. Theorem $(2.5)$], it then follows that there exists a neighborhood
of zero in $\Delta^n$ over which we may write
$$
	F(s) = e^{\Gamma(s)}.F
$$
relative to a unique $q_F$-valued holomorphic function $\Gamma(s)$ vanishing 
at the origin.  

\par To continue our calculations, we recall \cite{higgs} that for 
sufficiently small values of $s$, one may write 
$$
	I^{p,q}_{(F(s),W)} = e^{\Gamma(s)}e^{-\phi(s)}.I^{p,q}_{(F,W)}
$$
relative to a $C^{\infty}$ function $\phi(s)$ which has a first order
Taylor series expansion given by the the formula
$$
	\phi(s) = -\pi(\overline{L(s)})
$$
where $\pi:Lie(\Gc)\to Lie(\Gc^F)$ denotes projection with respect to the
the decomposition 
$$
	Lie(\Gc) = Lie(\Gc^F)\oplus q_F
$$
and 
$$
	L(s) = \sum_{j=1}^n\, \frac{\pd\Gamma(0)}{\pd s_j} s_j
$$
denotes the linearization of $\Gamma(s)$ about $s=0$.

\par To determine how $\nabla$ acts upon a smooth section $\sigma$ of 
$I^{p,q}_{(F(s),W)}$ observe that by virtue of the preceding remarks,
we may write
$$
	\sigma(s) = e^{\Gamma(s)}e^{-\phi(s)}\tilde\sigma(s)
$$
relative to a unique function $\tilde\sigma(s)$ taking values in the fixed
vector space $I^{p,q}_{(F,W)}$.  Consequently,
\begin{eqnarray*}
	\left.\nabla\sigma\right|_0 
	&=& \left.d e^{\Gamma(s)}e^{-\phi(s)}\right|_0\sigma(0) 
	     + \left.d\tilde\sigma\right|_0			\\
	&=& \left.dL\right|_0\sigma(0) +
	   \left.d\pi(\bar L)\right|_0\sigma(0) 
	 + \left.d\tilde\sigma\right|_0 
\end{eqnarray*}
In particular, since
$$
	L(s):\Delta^n\to\bigoplus_{s\leq 1}\, \gg^{-1,s}
$$
on account of the horizontality of $F(s)$,
$$
	\left.\nabla^{1,0}\sigma\right|_0 \in I^{p,q}_{(F,W)}
	+ \sum_{\ell\geq 0}\, I^{p-1,q+1-\ell}_{(F,W)},\qquad
	\left.\nabla^{0,1}\sigma\right|_0 \in I^{p+1,q-1}_{(F,W)}
	+ \sum_{\ell\geq 0}\, I^{p,q-\ell}_{(F,W)}
$$
\end{proof}

\par In order to state the next result, we note there is a natural functor
from the category of $\C$VMHS to the category of flat, filtered bundles 
which operates by replacing the bigrading $E=\oplus_{p,q}\, E^{p,q}$ of 
a $\C$VMHS by the associated filtration 
$\W_k = \bigoplus_{p+q\leq k}\, E^{p,q}$.

\setcounter{theorem}{20}					
\begin{theorem} The mixed Hodge metric $h$ carried by a complex variation of
graded-polarized mixed Hodge structure $E$ is filtered harmonic.
\end{theorem}
\begin{proof} By definition, a complex variation of graded-polarized mixed
Hodge structure consists of a complex variation of mixed Hodge structure
together with a collection of parallel hermitian forms on $Gr^{\W}$ which 
polarize the induced complex variations.  Accordingly, the associated 
decomposition
$$
	d_h = \bar\theta_0 + \bar\pd_0 + \pd_0 + \theta_0
$$
given by equation $(3.8)$ has the additional property that 
$$
	\bar\theta_0:E^{p,k-p}\to E^{p+1,k-p-1}\otimes\E^{0,1},\qquad
	\bar\pd_0:\E^0(E^{p,k-p})\to\E^{0,1}(E^{p,k-p})
$$
$$
	\theta_0:E^{p,k-p}\to E^{p-1,k-p+1}\otimes\E^{1,0},\qquad
	    \pd_0:\E^0(E^{p,k-p})\to\E^{1,0}(E^{p,k-p})
$$
On the other hand, by condition $(b)$ of Definition $(3.17)$, if 
$$
	\nabla = d_h + \tau_- + \theta_-
$$
denotes the decomposition of $\nabla$ defined by equation $(3.6)$, then
$$
	\tau_-:E^p\to E^p\otimes\E^{0,1},\qquad
	\theta_-:E^p\to E^{p-1}\otimes\E^{1,0}
$$
and hence the decomposition
$$
	\nabla = \bar\theta_0 + \bar\pd + \pd_0 + \theta
$$
obtained by setting $\bar\pd = \pd_0 + \tau_-$ and 
$\theta = \theta_0 + \theta_-$ coincides with usual decomposition of $E$ 
defined by equation $(3.12)$ via the system of Hodge bundles 
$E^p=\oplus_q\, E^{p,q}$.  In particular, $\bar\pd + \theta$ is an operator 
of Higgs--type.
\end{proof}

\par Moreover, in analogy with \cite{Simpson}, we have the following 
correspondence between filtered harmonic metrics and complex variations
of graded-polarized mixed Hodge structure:

\setcounter{theorem}{21}				    
\begin{theorem} Let $X$ be a compact K\"ahler manifold and $E\to X$ a flat,
filtered bundle.  Then, a filtered harmonic metric $h$ on $E$ underlies a
complex variation of graded-polarized mixed Hodge structure if and only if
\begin{description}
\item{$(a)$} $(E,\bar\pd+\theta)$ is a fixed point of the $\C^*$-action
$$
	  (E,\bar\pd+\theta)\mapsto(E,\bar\pd + \l\theta)	
$$
\item{$(b)$} The resulting decomposition of $E$ into a system of Hodge bundles
$$
	E = \bigoplus_p\, E^p				
$$
is preserved by the grading $\Y_h$.
\end{description}
\end{theorem}
\begin{proof} According to Theorem $(3.21)$, the mixed Hodge metric of a
complex variation of graded-polarized mixed Hodge structure is filtered
harmonic.  Conversely, given a filtered harmonic metric for which the above two
conditions hold one defines
\setcounter{equation}{22}
\begin{equation}
	E^{p,k-p} = E_k(\Y_h)\cap E^p				
\end{equation}

\par To prove the bigrading $(3.23)$ does indeed define a complex variation of 
graded-polarized mixed Hodge structure for which the associated mixed Hodge 
metric equals the given hermitian metric $h$, we note that for each index $k$,
the decomposition 
\begin{equation}
	E_k(\Y_h) = \bigoplus_p\, E^p_k,\qquad E^p_k = E^{p,k-p}  
\end{equation}
defines a system of Hodge bundles with respect to $\bar\pd_0 + \theta_0$.  
Indeed, this is automatic given conditions $(a)$ and $(b)$ since 
$\bar\pd_0 + \theta_0$ is just the component of $\bar\pd + \theta$ which 
preserves $\Y_h$.  On the other hand, by definition, $\bar\pd_0 + \theta_0$ 
coincides with the Higgs bundle structure obtained by applying Construction 
$(3.2)$ to the pair $(h,d_h)$.  Consequently, by Corollary $(3.16)$, there 
exists a unique complex variation of polarized Hodge structure on $E_k(\Y_h)$ 
for which $h$ is the Hodge metric, $d_h$ is the flat connection and equation 
$(3.24)$ represents the decomposition of $E_k(\Y_h)$ into a system of Hodge 
bundles.  Therefore, the decomposition
$$
	d_h = \bar\theta_0 + \bar\pd_0 + \pd_0 + \theta_0
$$
given by equation $(3.8)$ has the additional property that
$$
	\bar\theta_0:E^{p,k-p}\to E^{p+1,k-p-1}\otimes\E^{0,1},\qquad
        \bar\pd_0:\E^0(E^{p,k-p})\to\E^{0,1}(E^{p,k-p})
$$
$$
	\theta_0:E^{p,k-p}\to E^{p-1,k-p+1}\otimes\E^{1,0},\qquad
        \pd_0:\E^0(E^{p,k-p})\to\E^{1,0}(E^{p,k-p})
$$
Returning to condition $(a)$, it then follows that if 
$$
	\nabla = d_h + \tau_- + \theta_-
$$
denotes the decomposition of $\nabla$ defined by equation $(3.6)$ then 
$$
	\tau_-:E^p\to E^p\otimes\E^{0,1},\qquad
	\theta_-:E^p\to E^{p-1}\otimes\E^{1,0}
$$
since $\tau_- = \bar\pd - \bar\pd_0$ and $\theta_- = \theta-\theta_0$ on 
account of the filtered harmonicity of $h$.  Consequently, 
$$
	\nabla:\E^0(E^{p,k-p})\to\E^{0,1}(E^{p+1,k-p-1})\oplus
	\E^{0,1}(E^p)\oplus\E^{1,0}(E^{p,k-p})\oplus\E^{1,0}(E^{p-1})
$$
\end{proof}


\setcounter{theorem}{0}
\setcounter{lemma}{0}
\setcounter{definition}{0}
\setcounter{corollary}{0}
\setcounter{equation}{0}
\setcounter{example}{0}

\section{Admissibility Criteria}

\par Let $\Cal V\to S$ be a variation of graded-polarized mixed
Hodge structure and $\bar S$ be a compactification for which the
divisor $D=\bar S-S$ has at worst normal crossings.  Then, in contrast 
to the pure case, the associated period map $\varphi:S\to\M/\Gamma$ 
may have irregular singularities along $D$, as can be seen by 
considering the simplest of Hodge--Tate variations (cf. \S 2).
\vskip 10pt

\par To rectify this problem, let us return to our prototypical 
example
\begin{equation}
	\Cal V = R^k_{f\ast}(\C)		 
\end{equation}
defined by a family of algebraic varieties $f:X\to\Delta^*$ with 
unipotent monodromy.  Then, as discussed in \cite{SZ}, the limiting 
Hodge filtration $F_{\infty}$  of $\V$ exists.

\par Thus, comparing the asymptotic behavior of $(4.1)$ and $(4.2)$, 
we see that a minimal condition required in order for an abstract 
variation $\Cal V\to\Delta^*$ with unipotent monodromy to be akin to 
a geometric variation is that:
\begin{description}
\item{(i)}  \hph{a}The limiting Hodge filtration of $\V$ exists (in the
sense of \S 2).
\end{description}
Moreover, as shown by Deligne using $\ell$-adic techniques \cite{Deligne3}, 
the geometric variations $(4.2)$ are subject to a subtle condition which 
greatly restricts their local monodromy.  Namely, if $W$ denotes the constant
value of the weight filtration of $\V$, then: 
\begin{description}
\item{(ii)} \hph{a}The relative weight filtration $\rel W = \rel W(N,W)$
of the monodromy logarithm $N$ of $\V$ exists.
\end{description}
Consequently, on the basis of such considerations, we shall adopt the
terminology of \cite{SZ} and call an abstract variation $\Cal V\to\Delta^*$
with unipotent monodromy {\it admissible} provided it satisfies conditions 
$(i)$ and $(ii)$. 

\begin{nonumtheorem}[Deligne, \cite{SZ}] The limiting Hodge filtration of 
an admissible variation of graded-polarized mixed Hodge structure 
$\V\to\Delta^*$ with unipotent monodromy pairs with the relative weight 
filtration of $\V$ to define a mixed Hodge structure for which $N$ is 
morphism of type $(-1,-1)$.
\end{nonumtheorem}

\par To state main result of this section, let $\rel Y$ be a grading of the 
relative weight filtration $\rel W = \rel W(N,W)$ which is {\it compatible} 
with $N$ and $W$ in the following sense:
\begin{itemize}
\item $\rel Y$ preserves $W$.
\item $[\rel Y,N] = -2N$.
\end{itemize}
Then, given any grading $Y$ of $W$ preserving $\rel W$, we may construct
an associate $sl_2$ representation 
\begin{equation}
	(N_0,\rel Y-Y,N_0^+)				
\end{equation}
on the underlying vector space of $\rel W$ by decomposing $N$ as
\begin{equation}
	N = N_0 + N_{-1} + N_{-2} + \cdots		
\end{equation}
according to the eigenvalues of $ad\, Y$ and applying Theorem $(2.15)$
to the pair 
$$
	(N_0,\rel Y-Y)
$$

\setcounter{theorem}{4} 			
\begin{theorem}[Deligne, \cite{Deligne}] Let $\rel Y$ be a grading of 
$\rel W(N,W)$ which is compatible with $N$ and $W$.  Then, there exists a 
unique grading $Y$ of $W$ such that $Y$ preserves $\rel W$ and
\setcounter{equation}{5}
\begin{equation}
	[N-N_0,N_0^+] = 0				
\end{equation}
\end{theorem}
\begin{proof} We begin by selecting a grading $Y_0$ of $W$ which 
preserves $\rel W$, and recalling that by \cite{CKS}, the group
$$
	G^{\circ}
	 = \{\, g\in GL(V)^{\rel Y} \mid (g-1)(W_k)\subseteq W_{k-1} \,\}  
$$
acts simply transitively on the set of all such gradings $\Y(\rel Y,W)$.  
Next, to compute how this transitive action changes the associated $sl_2$ 
representations $(4.3)$, write
$$
	N=N_0 + N_{-1} + N_{-2} + \cdots
$$
relative to $Y\in\Cal Y(\rel Y,W)$, and observe that
\begin{equation}
	g\in G^{\circ} \implies Gr(Ad(g)\,N_0) = Gr(N)	
\end{equation}
Consequently, by $(4.7)$, the decomposition
$$
	N = N_0^* + N_{-1}^* + N_{-2}^* + \cdots
$$
of $N$ relative to $g.Y$ must have $N_0^*=Ad(g)\,N_0$, hence
\begin{equation}
   (N_0,\rel Y-Y,N_0^+) \stackrel{g}\rightarrow Ad(g)(N_0,\rel Y-Y,N_0^+)
\end{equation}

\par To prove the existence of a grading $Y\in\Cal Y(\rel Y,W)$ which satisfies
$(4.6)$, we proceed by induction and construct a sequence of gradings 
$Y_0,Y_1,\dots,Y_n$ in $\Cal Y(\rel Y,W)$ terminating at $Y$, 
by the requirement that
\begin{equation}
	[N-N_0,N_0^+] :W_j\to W_{j-\ell-1}	
\end{equation}
upon decomposing $N$ relative to $ad(Y_{\ell})$.  To check the validity of 
this algorithm, let us suppose the gradings $Y_0,Y_1,\dots,Y_{k-1}$ have 
been constructed, and write 
$$
	Y_k = g.Y_{k-1},\qquad g\in G^{\circ}	
$$
Then, by virtue of equation $(4.8)$, we can reformulate condition $(4.9)$ as
the requirement that 
\begin{equation}
	[Ad(g^{-1})N-N_0,N_0^+]:W_j \to W_{j-k-1} 	
\end{equation}
with 
$$
	N = N_0 + N_{-1} + N_{-2} + \cdots
$$
denoting the decomposition of $N$ relative to $ad(Y_{k-1})$.

\par Now, by virtue of the fact that
$$
	[N-N_0,N_0^+]:W_j\to W_{j-k}
$$
(upon decomposing $N$ relative to $ad(Y_{k-1})$), it is natural to assume
that our element $g$ is of the form 
$$
	g = 1 + \g_{-k}+\g_{-k-1} + \cdots
$$
(again, relative to the eigenvalues of $ad(Y_{k-1})$).  Imposing condition 
$(4.10)$, it follows that we may take $g$ to be of this form if and only 
there exists a solution $\g_{-k}$ of the equation:
\begin{equation}
	[N_{-k}+[N_0,\g_{-k}],N_0^+] = 0		
\end{equation}
To find such an element $\g_{-k}$, one simply observes that
$$
	End(V) = \Im(ad\,N_0)\oplus\ker(ad N_0^+)
$$

\par To prove the grading $Y$ so constructed is unique, suppose that $Y'$ is 
another such grading and write $Y' = g.Y$ with
$$
	g = 1 + \g_{-k}+\g_{-k-1} + \cdots\in G^{\circ},\qquad \g_{-k}\neq 0
$$
relative to $ad(Y)$.  Therefore, upon applying $Ad(g^{-1})$ to both sides
of equation $(4.8)$, we see that
$$
	[Ad(g^{-1})N-N_0,N_0^+] = 0
$$
and hence, we also must have 
\begin{equation}
	[N_{-k}+[N_0,\g_{-k}],N_0^+] = [[N_0,\g_{-k}],N_0^+] = 0	
\end{equation}
since $g=1+\g_{-k}+\g_{-k-1}+\cdots$.

\par In addition, by combining the observation that $\rel Y$ preserves the
eigenspaces  of $Y$ with the fact that $g\in G^{\circ}$ fixes $\rel Y$,
it follows that we must also have
$$
	[\g_{-k},\rel Y] = 0 
$$
Thus, $\g_{-k}$ is a solution to $(4.12)$ which is of weight $k>0$ relative
to the representation $(ad(N_0),ad(\rel Y-Y),ad(N_0^+))$.  By standard 
$sl_2$ theory, we therefore have $\g_{-k} = 0$, contradicting our 
our assumption that $g\neq 1$.
\end{proof}

\begin{remark} Deligne \cite{Deligne}: 
\begin{itemize}
\item Since $[N-N_0,N_0^+]=0$, each non-zero term $N_{-k}$ with $k>0$ 
      appearing in the decomposition $(4.4)$ defines an irreducible 
      representation of $sl_2(\C)$ of highest weight $k-2$ via the 
      adjoint action of $(N_0,\rel Y-Y,N_0^+)$.  In particular, 
      $N_{-1} = 0$.
\item The construction of Theorem $(4.5)$ is compatible with tensor products, 
      directs sums and duals.
\item The construction is functorial with respect to morphisms.\newline
{\it Sketch of Proof}.  {\rm Let $f=\sum_{k\geq 0}\,f_{-k}$ be 
the decomposition of our morphism $f:V_1\to V_2$ relative 
to the induced grading of $Hom(V_1,V_2)$.  Then, in analogy 
with the proof of Theorem $(4.5)$, direct computation shows 
that relative to the induced representation $(N_0,H,N_0^+)$ 
on $Hom(V_1,V_2)$ we must have both  
$$
	[N_0^+,[N_0,f_{-k}]]=0
$$ 
and $H(f_{-k}) = kf_{-k}$.}
\end{itemize}
\end{remark}

\par To relate this construction of Deligne to admissible variations, let us 
suppose we have such a variation $\Cal V\to\Delta^*$ with limiting mixed 
Hodge structure $(F_{\infty},\rel W)$.  Then, because $N$ is a morphism of 
$(F_{\infty},\rel W)$ of type $(-1,-1)$, we may apply Theorem $(4.5)$ with 
\begin{equation}
		\rel Y = Y_{(F_{\infty},\rel W)}   
\end{equation}
to obtain an associate grading:
$$
	Y = Y(F_{\infty},W,N)
$$

\par Thus, motivated by these observations, let us call a 
triple $(F,W,N)$ {\it admissible} provided:
\begin{description}
\item{(a)} \hph{a}The relative weight filtration $\rel W = \rel W(N,W)$ exists.
\item{(b)} \hph{a}The pair $(F,\rel W)$ induces a mixed Hodge structure on each
	   $W_k$,  relative to which $N$ is $(-1,-1)$-morphism.
\end{description}
Then, mutatis mutandis, we obtain an associate grading 
\begin{equation}
	Y=Y(F,W,N)					
\end{equation}
for each admissible triple $(F,W,N)$.

\par To study the dependence of the grading $(4.15)$ upon the filtration $F$,
recall that via the splitting operation described in \cite{CKS}, we can pass 
from an arbitrary mixed Hodge structure $(F,\rel W)$ to a mixed Hodge 
structure $(\hat F,\rel W) = (e^{-i\delta}.F,\rel W)$ split over $\R$.  
Moreover, a moment of thought shows 
$$
	e^{-i\delta}:(\rel Y,F,\rel W)\to (\rel\hat Y,N,W)
$$
to be a morphism of admissible triples, hence
\begin{equation}
	Y(F,W,N) = e^{i\delta}.Y(\hat F,W,N)	        
\end{equation}
by the functoriality of Theorem $(4.5)$.

\par To understand the power of this reduction to the split over $\R$ case,
observe that every mixed Hodge structure which splits over $\R$ and has 
$N$ as a $(-1,-1)$-morphism can be built up from the following two sub-cases:
\begin{description}
\item[$(1)$] $V=\bigoplus\, I^{p,p}$ and $N(I^{p,p})\subseteq I^{p-1,p-1}$.
\item[$(2)$] $V=I^{p,0}\bigoplus I^{0,p}$ and $N=0$.
\end{description}
via the operations of tensor product and direct sums.  In particular, 
because Theorem $(4.5)$ is compatible with these operations, we obtain
the following result:

\setcounter{theorem}{16}			   
\begin{theorem} Let $(F,W,N)$ be an admissible triple. Then, the grading 
$Y=Y(F,W,N)$ preserves the filtration $F$.
\end{theorem}
\begin{proof} By virtue of the previous remarks, it will suffice to check the
two sub-cases enumerated above.  To verify the assertion for case (1), 
note that $I^{p,p}$ is the weight 2p eigenspace of $\rel Y$, and recall that
by definition we must have
\begin{equation}
	F^p = \bigoplus_{a\geq p}\, I^{p,p}		
\end{equation}
Consequently, the commutativity relation $[\rel Y,Y] = 0$ implies that 
$$
	Y(I^{p,p})\subseteq I^{p,p}			
$$
and hence $Y$ preserves $F$ by equation $(4.18)$.

\par Regarding the second case, observe that because $N=0$, $W$ must be the
trivial filtration
$$
	0=W_{2p-1}\subset W_{2p} = V
$$
in order for the relative weight filtration $\rel W = \rel W(N,W)$ to exist.
Thus, $Y$ must be the trivial grading on $V$ of weight $2p$.
\end{proof}

\par Returning now to the work of Schmid, let us recall the following basic
result concerning the monodromy weight filtration \cite{Schmid}:

\setcounter{lemma}{18}					
\begin{lemma} Let $V$ be a finite dimensional $\C$-vector space endowed
with a non-degenerate bilinear form $Q$, and suppose $N$ is a nilpotent
endomorphism of $V$ which acts by infinitesimal isometries of $Q$.
Then, any semisimple endomorphism $H$ of $V$ which satisfies $[H,N] = -2N$
is also an infinitesimal isometry of $Q$.
\end{lemma}
\begin{proof} There are two key steps:
\begin{description}
\item{$(1)$} \hph{a} Show the monodromy weight filtration
	\setcounter{equation}{19}
	\begin{equation}
	    0\neq W_{-\ell}(N)\subseteq\dots W_{\ell}(N) = V	
	\end{equation}
	is self-dual with respect to $Q$, i.e. 
	$W_j(N) = W_{-j-1}(N)^{\perp}$.
\item{$(2)$} \hph{a} Via semisimplicity, assume the pair $(N,H)$ defines an 
	$sl_2$ representation $\{e_k,e_{k-2},\dots,e_{2-k},e_{-k}\}$ of 
	highest weight $k$.  Imposing self-duality, it then follows that 
	$Q(e_i,e_j) = 0$ unless $i+j=0$, hence $H$ is an infinitesimal
	isometry of $Q$.
\end{description}
\end{proof}

\setcounter{corollary}{19}	          
\begin{corollary} Let $\Cal V\to\Delta^*$ be an admissible variation with
limiting data $(F_{\infty},W,N)$ and graded-polarizations $\{\Cal S_k\}$.
Then, the semi-simple element $\rel Y-Y$ constructed by Theorem $(4.5)$ acts
on $Gr^W$ by infinitesimal isometries.
\end{corollary}

\setcounter{theorem}{0}
\setcounter{lemma}{0}
\setcounter{definition}{0}
\setcounter{corollary}{0}
\setcounter{equation}{0}
\setcounter{example}{0}

\section{The Nilpotent Orbit Theorem}

\par In this section, we state and prove a precise analog of Schmid's 
Nilpotent Orbit Theorem for admissible variations of graded-polarized
mixed Hodge structure $\V\to\Delta^*$ with unipotent monodromy.
\vskip 10pt

\par To this end, we recall from \S 2 that the period map 
$$
	\varphi:\Delta^*\to\M/\Gamma
$$ 
of such a variation lifts to a holomorphic, horizontal map 
$F(z):U\to\M$ which makes the following diagram commute:
$$
\begin{CD}
                  U        @> F >>    \M        \\
 	@Vs=\exp(2\pi i z)VV                    @VVV     \\
                  \Delta^* @> \phi >> \M/\Gamma
\end{CD}
$$
In particular, as noted in \S 2, the map $F(z)$ then descends to an 
\lq\lq untwisted period map\rq\rq{}
$$
	\psi(s):\Delta^*\to\check\M
$$
on account of the quasi-periodicity condition $F(z+1) = e^N.F(z)$.

\begin{theorem}{\bf [Nilpotent Orbit Theorem]} Let $\psi(s)$ be the untwisted
period map of an admissible variation of graded-polarized mixed Hodge
structure $\V\to\Delta^*$ with unipotent monodromy.  Then,
\begin{description}
\item{(i)}   \hph{a}The limiting Hodge filtration 
	     $F_{\infty} = \lim_{s\to 0}\,\psi(s)$ of $\V$ exists, and is
	     an element of $\check\M$.
\item{(ii)}  The nilpotent orbit $F_{nilp}(z) = \exp(zN).F_{\infty}$
             extends to a holomorphic, horizontal map $\C\to\check\M$.
             Moreover, there exists $\a>0$ such that $F_{nilp}(z)\in\M$
             whenever $\Im(z)>\a$.
\item{(iii)} Let $d_{\M}$ denote the distance function determined by the
	     mixed Hodge metric $h$.  Then, there exists constants $K$ 
	     and $\b$ such that:
             $$
                 d_{\M}(F(z),F_{nilp}(z))
                 \leq K\Im(z)^{\b}\exp(-2\pi\Im(z))
             $$
	     for all $z\in U$ with $\Im(z)$ sufficiently large.
\end{description}
\end{theorem}

\par Regarding the proof of Theorem $(5.1)$, observe that part $(i)$ is a 
direct consequence of the admissibility of $\Cal V$ (cf. \S 2).  Likewise, 
part $(ii)$ is a direct consequence of Schmid's Nilpotent Orbit Theorem, 
applied to the variations of pure, polarized Hodge structure carried by 
$Gr^W$.  

\par To prove part $(iii)$, recall from \S 2 that $\check\M$ is a complex
manifold upon which the Lie group
$$
	\Gc = \{\, g\in GL(V)^W \mid Gr(g)\in Aut(\Cal S,\C)\, \}
$$
acts transitively.  Therefore, given any element $F\in\check\M$ and a vector 
space decomposition
\setcounter{equation}{1}
\begin{equation}
        Lie(\Gc) = Lie(\Gc^F)\oplus q			
\end{equation}
the map $u\in q\mapsto e^u.F$ will be a biholomorphism from a neighborhood 
of zero in $q$ onto a neighborhood of $F$ in $\check\M$.  In particular, 
upon setting $F=F_{\infty}$, and shrinking $\Delta$ as necessary, we see that
each choice of decomposition $(5.2)$ determines a corresponding holomorphic map
$$
       \Gamma(s):\Delta\to q,\qquad \Gamma(0) = 0
$$
via the rule:
\begin{equation}
        e^{\Gamma(s)}.F_{\infty} = \psi(s)		
\end{equation}

\par To make a choice of decomposition $(5.2)$, we shall follow the methods of
\cite{higgs} and use the limiting mixed Hodge structure of $\Cal V$ to define 
a grading of $Lie(\Gc)$:

\setcounter{lemma}{3} 					
\begin{lemma} The limiting mixed Hodge structure $(F_{\infty},\rel W)$ of an
admissible variation of graded-polarized mixed Hodge structure defines a 
canonical, graded, nilpotent Lie algebra
\setcounter{equation}{4}
\begin{equation}
	q_{\infty} = \bigoplus_{a<0}\,\wp_a		
\end{equation}
with the following additional property:  As complex vector spaces,
\begin{equation}
	Lie(\Gc) = Lie(\Gc^{F_{\infty}})\oplus q_{\infty}	
\end{equation}
\end{lemma}
\begin{proof} The details may be found in \S 6 of \cite{higgs}.  However, 
the idea of the proof is relatively simple:  An element $\a\in Lie(\Gc)$
belongs to the subspace $\wp_a$ if and only if 
$$
	\a:I^{p,q}_{(F_{\infty},\rel W)} 
	\to \bigoplus_b\, I^{p+a,b}_{(F_{\infty},\rel W)}
$$
To verify the decomposition $(5.6)$, we observe that 
\begin{description}
\item{$(1)$} As a vector space, $Lie(\Gc) = \oplus_a\,\wp_a$.
\item{$(2)$} By definition, 
	     $F^p_{\infty} = \oplus_{a\geq p}\, I^{a,b}_{(F_{\infty},\rel W)}$
	     and hence
	     $Lie(\Gc^{F_{\infty}}) = \bigoplus_{a\geq 0}\,\wp$.
\end{description}
\end{proof}
	
\par To establish the distance estimate $(iii)$, let us now record the 
following result:

\setcounter{lemma}{6}				    
\begin{lemma} Let $Y$ be a grading of $W$ which is defined over $\R$ and
$y$ be a positive real number.  Then, 
$$
	\a\in\R \implies y^{\a Y}:\M\to\M
$$
Moreover, if $\a<0$, $d_{\M}$ denotes the Riemann distance on $\M$ defined by 
the mixed Hodge metric (cf. \S 2) and $L$ denotes the length of $W$ then:
\setcounter{equation}{7}
\begin{equation}
	F_1,F_2\in\M\implies d_{\M}(y^{\a Y}.F_1,y^{\a Y}.F_2) 
	\leq y^{-\a(L-1)}d_{\M}(F_1,F_2)		
\end{equation}
\end{lemma}	  
\begin{proof} To prove that $y^{\a Y}:\M\to\M$, one simply checks that 
$y^{\a Y}$ acts as scalar multiplication by $y^{\a k}$ on $Gr^W_k$, and
hence $(y^{\a Y}.F)Gr^W_k = FGr^W_k$.

\par To verify equation $(5.8)$, we note that if $Y_F$, $F\in\M$ denotes 
the grading of $W$ which acts as multiplication by $(p+q)$ on 
$I^{p,q}_{(F,W)}$ then 
$$
	Y_{y^{\a Y}.F} = Ad(y^{\a Y}) Y_F
$$ 
since $Y$ is defined over $\R$, and hence
$$
	v\in E_k(Y_F)\implies 
	||y^{\a Y}v||_{y^{\a Y}.F} = y^{\a k}||v||_F
$$
Upon transferring these computations to the induced metrics on $Lie(\Gc)$, it
then follows that
$$
	T\in E_{-\ell}(ad\, Y_F) \implies 
	||Ad(y^{\a Y})T||_{y^{\a Y}.F} = y^{-\a\ell}||T||_F
$$
In particular, $||Ad(y^{\a Y})T||_{y^{\a Y}.F} \leq y^{-\a(L-1)}||T||_F$
since $\a<0$, and hence  
$$
	d_{\M}(y^{\a Y}.F_1,y^{\a Y}.F_2)\leq y^{-\a(L-1)}d_{\M}(F_1,F_2)
$$
\end{proof}

\begin{remark} In connection with the proof of Theorem $(5.9)$, we note the
following result: Let $(M,g)$ be a Riemannian manifold which is an open subset
of a manifold $\check M$ upon which a Lie group $\Cal G$ acts transitively, 
and $|\cdot|$ be a norm on $Lie(\Cal G)$.  Then given any point $F_0\in M$ 
there exists a neighborhood $S$ of $F_0$ in $M$, a neighborhood $U$ of zero 
in $Lie(\Cal G)$ and a constant $K>0$ such that 
$$
        u\in U, F\in S\implies e^u.F\in M \quad\text{and}\quad
        d_M(e^u.F,F) < K |u|
$$
\end{remark}

\setcounter{theorem}{8} 			      
\begin{theorem}{\bf [Distance Estimate]} Let $F(z):U\to\M$ be a lifting of 
the period map $\varphi$ of an admissible variation of graded-polarized 
mixed Hodge structure $\V\to\Delta^*$ with unipotent monodromy logarithm
$N$ and limiting mixed Hodge structure $(F_{\infty},\rel W)$.  Then, given 
any $G_{\R}$ invariant metric on $\M$ which obeys $(5.8)$, there exists 
constants $K$ and $\b$ such that 
$$
        Im(z)>>0 \implies d_{\M}(F(z),e^{zN}.F_{\infty}) < Ky^{\b} e^{-2\pi y}
$$
\end{theorem} 
\begin{proof} For simplicity of exposition, we shall first prove the result
under the additional assumption that our limiting mixed Hodge structure 
$(F_{\infty},\rel W)$ is split over $\R$.   Having made this assumption,
it then follows from the work of \S 4 that:
\begin{description}
\item{$(1)$} The associate gradings $\rel Y=Y_{(F_{\infty},\rel W)}$ and 
		$Y = Y(F_{\infty},W,N)$ are defined over $\R$.
\item{$(2)$} The endomorphism $\rel Y-Y$ is an element $Lie(G_{\R})$.
\item{$(3)$} The filtration $F_0 = e^{iN}.F_{\infty}$ is an element of 
	     $\M$. [This assertion is a consequence of Schmid's $SL_2$ Orbit 
	     Theorem, see \cite{Schmid} for details.]
\end{description}

\par Now, as discussed in \S 4, the fact that $N$ is a $(-1,-1)$-morphism
of the limiting mixed Hodge structure $(F_{\infty},\rel W)$ implies that
$$
	[\rel Y,N] = -2N
$$
and hence
$$
        e^{iyN} = y^{-\half\rel Y}e^{iN}y^{\half\rel Y}
$$ 
Consequently, because $\rel Y$ preserves $F_{\infty}$, 
$$
	e^{iyN}.F_{\infty} = y^{-\half\rel Y}e^{iN}y^{\half\rel Y}
			   =  y^{-\half\rel Y}e^{iN}.F_{\infty} 
			   = y^{-\half\rel Y}.F_0
$$

\par Next, we note that by equation $(5.3)$ and Lemma $(5.4)$
$$
	F(z) = e^{zN}.\psi(s) = e^{zN}e^{\Gamma(s)}.F_{\infty}
$$
and hence
\begin{eqnarray*}
    d_{\M}(e^{zN}e^{\Gamma(s)}.F_{\infty},e^{zN}.F_{\infty}) 
    &=& d_{\M}(e^{iyN}e^{\Gamma(s)}.F_{\infty},e^{iyN}.F_{\infty}) \\   
    &=& d_{\M}(y^{-\half\rel Y}e^{iN}y^{\half\rel Y}e^{\Gamma(s)}.F_{\infty},
           y^{-\half\rel Y}e^{iN}y^{\half\rel Y}.F_{\infty})
\end{eqnarray*}
In particular, upon setting
$$
        e^{\tilde\Gamma(z)} = Ad(e^{iN})Ad(y^{\half\rel Y})e^{\Gamma(s)}
$$ 
and recalling that $\rel Y$ preserves $F_{\infty}$, it follows that
\setcounter{equation}{9}
\begin{equation}
    d_{\M}(e^{zN}e^{\Gamma(s)}.F_{\infty},e^{zN}.F_{\infty}) 
    = d_{\M}(y^{-\half\rel Y}e^{\tilde\Gamma(z)}.F_0, y^{-\half\rel Y}.F_0)
\end{equation}

\par In addition, because $[\rel Y,Y]= 0$, 
$$
	y^{-\half \rel Y} = y^{-\half(\rel Y-Y)}y^{-\half Y}
$$
and hence
\begin{equation}  
   d_{\M}(y^{-\half\rel Y}e^{\tilde\Gamma(z)}.F_0, y^{-\half\rel Y}.F_0)
   = d_{\M}(y^{-\half Y}e^{\tilde\Gamma(z)}.F_0, y^{-\half Y}.F_0)
\end{equation}
since $\rel Y-Y$ is an element of $Lie(G_{\R})$.  Therefore, by equation
$(5.8)$, 
$$
        d_{\M}(e^{zN}e^{\Gamma(s)}.F_{\infty},e^{zN}.F_{\infty})
        \leq y^{\half(L-1)}d_{\M}(e^{\tilde\Gamma(z)}.F_0,F_0)
$$
Consequently, by the remark which follows Lemma $(5.7)$, we have:  
$$
        d_{\M}(e^{zN}e^{\Gamma(s)}.F_{\infty},e^{zN}.F_{\infty})
        < K|\tilde\Gamma(z)|\text{Im}(z)^{\half(L-1)}
$$
for some $K>0$ and norm any fixed norm $|\cdot|$ on $gl(V)$.

\par To further analyze $\tilde\Gamma(z)$, decompose $\tilde\Gamma(z)$
according to the eigenvalues of $ad\,\rel Y$
$$
        \Gamma = \sum_{\ell}\,\Gamma_{[\ell]},
        \qquad [\rel Y,\Gamma_{[\ell]}] = \ell\Gamma_{[\ell]}
$$
Then
\begin{eqnarray*}
  \tilde\Gamma(z) &=& e^{i\,ad\,N}y^{\half ad\,\rel Y}\Gamma(s) 
   = e^{i\,ad\,N}y^{\half ad\,\rel Y}\sum_{\ell}\,\Gamma_{[\ell]}(s)  \\ 
   &=& \sum_{\ell}\,y^{\frac{\ell}{2}}e^{i\,ad\,N}\Gamma_{[\ell]}(s)
\end{eqnarray*}
since $y^{\half ad\rel Y}\Gamma_{[\ell]} =y^{\frac{\ell}{2}}\Gamma_{[\ell]}$. 
Therefore, denoting the maximal eigenvalue of $ad\,\rel Y$ on $q_{\infty}$
by $c$, we have
$$
        |\tilde\Gamma(z)| \leq O(y^{c/2}e^{-2\pi y})
$$
because $\Gamma(s)$ is a holomorphic function of $s=e^{2\pi iz}$ vanishing at
zero.  

\par In summary, we have proven that whenever the limiting mixed
Hodge structure $(F_{\infty},\rel W)$ of our admissible variation $\Cal V$
is split over $\R$ and $Im(z)$ is sufficiently large, the following distance
estimate holds:
$$
        d_{\M}(e^{zN}e^{\Gamma(s)}.F_{\infty},e^{zN}.F_{\infty})
        \leq K y^{\b}e^{-2\pi y},\qquad\b = \half(c+L-1)
$$

\par If the limiting mixed Hodge structure $(F_{\infty},\rel W)$ is not 
split over $\R$, we may obtain the same distance estimate by first applying 
Deligne's $\delta$ splitting:
$$
   F_{\infty} = e^{i\delta}.\hat F_{\infty}
$$ 
and then proceeding as above.  

\par More precisely, let $\rel Y,Y$ denote the gradings determined by 
$(\hat F_{\infty},\rel W)$ and $N$.  Then, equation $(5.10)$ becomes
$$
   d_{\M}(e^{zN}e^{\Gamma(s)}.F_{\infty},e^{zN}.F_{\infty}) = 
      d_{\M}(y^{-\half Y}e^{\tilde\Gamma(z)}.F_0(y), y^{-\half Y}.F_0(y))
$$
with
$$
        F_0(y) = e^{i\delta(y)}e^{iN}.\hat F_{\infty},\qquad
        \delta(y) = y^{\half ad\,\rel Y}\delta
$$
and 
$$
	\delta\in\lam_{(F_{\infty},\rel W)} \implies \delta(y) \sim O(y^{-1})
$$
Thus, applying equation $(5.8)$ and the remark the follows the proof of 
Lemma $(5.7)$, we have:
$$
        d_{\M}(e^{zN}e^{\Gamma(s)}.F_{\infty},e^{zN}.F_{\infty})
        \leq Ky^{\b}e^{-2\pi y}
$$
for $Im(z)$ sufficiently large.
\end{proof}

\par In regards to analogs of Schmid's $SL_2$ Orbit Theorem for admissible
variations $\V\to\Delta^*$ of graded-polarized mixed Hodge structure, the 
following result suggests that as soon as the weight filtration of $\V$ has 
length $L\geq 3$, the resulting nilpotent orbit $e^{zN}.F_{\infty}$ need
not admit an approximation by an auxiliary orbit $e^{zN}.\hat F_{\infty}$ 
for which the limiting mixed Hodge structure $(\hat F_{\infty},\rel W)$ is 
split over $\R$:

\setcounter{theorem}{10}			
\begin{theorem} Let $(F_{\infty},\rel W)$ be the limiting mixed Hodge 
structure associated to an admissible variation of graded-polarized mixed 
Hodge structure with unipotent monodromy and $(\hat F_{\infty},\rel W)$ be
a mixed Hodge structure of the form 
$$
	(\hat F_{\infty},\rel W) = (e^{-\sigma}.F_{\infty},\rel W),\qquad
	\sigma\in\ker(ad\,N)\cap\lam_{(F_{\infty},\rel W)}
$$
which is split over $\R$.  Then, there exist a positive constant $K$ such that
$$
	d_{\M}(e^{zN}.\hat F_{\infty},e^{zN}.F_{\infty}) 
	\leq K \Im(z)^{\frac{L-3}{2}}
$$
for all $z\in U$ with $\Im(z)$ sufficiently large.
\end{theorem}
\begin{proof} Let $\rel Y$ and $Y$ denote the gradings associated to the
split mixed Hodge structure $(\hat F_{\infty},\rel W)$ via the methods 
of \S 4.  Then, a quick review of the proof of Theorem $(5.9)$ shows that
upon setting
\begin{itemize}
\item $F_0 = e^{iN}.\hat F_{\infty}\in\M$.
\item $\sigma(y) = (y^{\half ad\,\rel Y})(\sigma)$.
\end{itemize}
we have both $e^{zN}.\hat F_{\infty} = e^{xN}y^{-\half\rel Y}.F_0$
and $e^{zN}.F_{\infty} = e^{xN}y^{-\half\rel Y}e^{\sigma(y)}.F_0$

\par In particular, since $\rel Y-Y\in Lie(G_{\R})$  and $Y$ is a grading of 
$W$ which is defined over $\R$:
\begin{eqnarray*}
	d_{\M}(e^{zN}.\hat F_{\infty},e^{zN}.F_{\infty})
	&=& d_{\M}( e^{xN}y^{-\half\rel Y}.F_0,
	       e^{xN}y^{-\half\rel Y}e^{\sigma(y)}.F_0)	\\
	&=& d_{\M}(y^{-\half(\rel Y-Y)}y^{-\half Y}.F_0,
	      y^{-\half(\rel Y-Y)}y^{-\half Y}e^{\sigma(y)}.F_0)	\\
	&=& d_{\M}(y^{-\half Y}.F_0, y^{-\half Y}e^{\sigma(y)}.F_0)	\\
	&\leq& y^{\half(L-1)}d_{\M}(F_0,e^{\sigma(y)}.F_0)
\end{eqnarray*}
Moreover,
$
	\sigma\in\lam_{(F_{\infty},\rel W)} \implies \sigma(y)\sim O(y^{-1})
$
and hence
$$
	d_{\M}(F_0,e^{\sigma(y)}.F_0) \leq K y^{-1}
$$  
Therefore,
$$
	d_{\M}(e^{zN}.\hat F_{\infty},e^{zN}.F_{\infty})
	\leq K \Im(z)^{\frac{L-3}{2}}	
$$
\end{proof}

\par To verify that the distance estimate of Theorem $(5.11)$ is sharp in
the case where $L=3$, let $\M\cong\C$ denote the classifying space of 
Hodge--Tate structures constructed in Example $(2.7)$, and $N$ be the 
nilpotent endomorphism of $V=\text{span}(e_2,e_0)$ defined by the rule 
$$
	N(e_2) = e_0,\qquad N(e_0) = 0
$$ 
Then, a short calculation shows that
\begin{itemize}
\item $\rel W = \rel W(N,W)$ exists and coincides with $W$.  
\item At each point $F\in\M$, $Lie(\Gc)=\gg^{-1,-1}_{(F,W)}
      =\text{span}_{\C}(N)$.  
\end{itemize}
In particular, given any point $F_{\infty}\in\M$, the map
$$
	z\mapsto e^{zN}.F_{\infty}
$$
is an admissible nilpotent orbit.  Moreover, since $\M$ is Hodge--Tate, the
resulting mixed Hodge metric on $\M$ is invariant under left translation by
$e^{\l N}$ for all $\l\in\C$.  Consequently, upon setting
$$
	F_{\infty} = e^{\sigma}.\hat F_{\infty}
$$
for some element $\sigma\in\lam_{(F,W)}$ and some point 
$\hat F_{\infty}\in\M_{\R}$, it then follows that
$$
	d_{\M}(e^{zN}.\hat F_{\infty},e^{zN}.F_{\infty})
	= d_{\M}(e^{zN}.\hat F_{\infty},e^{zN}e^{\sigma}.\hat F_{\infty}) 
	= d_{\M}(\hat F_{\infty},e^{\sigma}.\hat F_{\infty})
$$

\end{document}